\documentclass[a4paper,11pt]{article}
\usepackage{graphicx,soul,comment,booktabs} 
\usepackage[table]{xcolor}
\definecolor{dred}{rgb}{0.6,0,0}
\definecolor{dpurple}{HTML}{A020F0}
\definecolor{dblue}{rgb}{0,0,1}
\definecolor{hlcolor}{rgb}{1,1,0.8}
\definecolor{control}{HTML}{41AB5D}
\usepackage[top=3cm,bottom=3cm,left=3cm,right=3cm]{geometry}
\usepackage{amsmath,amssymb,amsthm,bm,framed,natbib,authblk}
\usepackage{lineno}
\usepackage[colorlinks=true,urlcolor=dblue, linkcolor=dblue,citecolor=dred]{hyperref}
\usepackage[nameinlink]{cleveref}
\usepackage[labelfont=bf,labelsep=period]{caption}
\usepackage{graphicx}

\usepackage{libertine}
\renewcommand\b\bm

\title{Automatic differentiation of Sylvester, Lyapunov, and algebraic Riccati equations}
\author[1]{\normalsize \bfseries Ta-Chu Kao${}^{\rm @}$} 
\author[1]{\bfseries Guillaume Hennequin} 

\renewcommand\bar\overline

\affil[1]{Computational and Biological Learning Lab, Department.~of Engineering, University~of~Cambridge, Cambridge, UK} 
\date{\normalsize ${}^{\rm @}$ tck29@cam.ac.uk\\ \vspace*{1cm} \today}
\sethlcolor{hlcolor}
\begin{document}
\maketitle

\begin{abstract}
  Sylvester, Lyapunov, and algebraic Riccati equations are the bread and butter of control theorists.
  They are used to compute infinite-horizon Gramians, solve optimal control problems in continuous or discrete time, and design observers.
  While popular numerical computing frameworks (e.g., \scalebox{0.85}{\texttt{scipy}}) provide efficient solvers for these equations, these solvers are still largely missing from most automatic differentiation libraries.
  Here, we derive the forward and reverse-mode derivatives of the solutions to all three types of equations, and showcase their application on an inverse control problem.
\end{abstract}

\section{Introduction}
\label{sec:introduction}
In recent years, automatic differentiation (AD) has become the backbone of most modern machine learning applications and an important tool for scientific enquiry \citep{Baydin2017}.
The success of AD owes in part to the myriad of forward and reverse-mode derivatives that have been derived for matrix/tensor operations including most known linear operators and factorizations \citep{Giles2008}.
These are commonly available in popular automatic differentiation libraries such as PyTorch \citep{Paszke2019}, Tensorflow \citep{Paszke2019}, Jax \citep{Bettencourt2019}, Diffsharp \citep{Baydin2017}, and Zygote \citep{Innes2018}.

In this technical note, we derive the forward and reverse-mode gradients for solutions to three important types of matrix equations used extensively in control theory: Sylvester, Lyapunov, and algebraic Riccati equations.
Although efficient solvers are readily available in popular scientific programming libraries, to the best of our knowledge they remain missing from most automatic differentiation packages (though CasADi seems to provide Lyapunov solvers; \citealp{Gillis2015, Andersson2019}).
We have added these solvers to Owl~\citep{Wang2017}, a numerical library written in OCaml with a full-featured AD module, and we hope that this technical note will enable rapid integration in other popular AD packages.

\section{Preliminaries}
\label{sec:preliminaries}
We use bold upper-case letters to denote matrices and bold lower-case letters to denote vectors.
We use $\b{X}^T$, $\b{X}^{-1}$ and $\text{tr}(\b{X})$ to denote the transpose, inverse, and trace of a matrix $\b{X} \in \mathbb{R}^{n \times n}$. 
Following \citet{Giles2008}, we use $\bar{\b{X}}$ to denote the adjoint of matrix $\b{X}$ w.r.t.\ a scalar output $\ell$, i.e.\ $\partial \ell / \partial \b{X}$.
We also use $\dot{\b{X}}$ to denote the tangent of a matrix $\b{X}$ w.r.t. some scalar input $s$, i.e.\ $\partial \b{X}/\partial s$.

Consider a differentiable solver $\b{P} = f(\b{A}, \b{B}, \ldots)$.
In forward mode, the input tangents $(\dot{\b{A}}, \dot{\b{B}}, \ldots)$ are known and used to calculate $\dot{\b{P}}$.
In reverse mode, the solution's adjoint $\bar{\b{P}}$ is known and used to update the input's adjoints $(\bar{\b{A}}, \bar{\b{B}}, \ldots)$.
This note consists mainly in deriving the rules by which these tangents and adjoints should be computed for various types of matrix equation solvers. 
All our results are summarized in \Cref{sm:table:forward,sm:table:reverse}.

\section{Sylvester equations}
\label{sec:sylvester}
\subsection{Continuous time Sylvester equation}

The continuous time Sylvester equation takes the form
\begin{equation} \label{eq:csylv}
  \b{A}\b{P}+ \b{P}\b{B} + \b{C} = \b{0},
\end{equation}
where $\b{A}$, $\b{P}$, $\b{B}$, and $\b{C}$ are all square matrices that have the same dimensions.
There is a unique solution $\b{P}$ that satisfies \Cref{eq:csylv}
for a given $\b{A}$, $\b{B}$, and $\b{C}$ if and only if $\b{A}$ and $\b{B}$ have no common eigenvalue.
As we shall see in the following, the matrix derivatives for the solution to the Sylvester equation is also unique under the same condition on the spectrum of $\b{A}$ and $\b{B}$.

\paragraph{Forward mode}
Taking the differential on both sides of \Cref{eq:csylv}, we have
\begin{equation}
  \b{A}\, \b{dP} + \b{dP}\, \b{B} + 
  \big ( \b{dA} \, \b{P} + \b{P} \, \b{dB} + \b{dC} \big) = \b{0}.
\end{equation}
The forward mode derivative $\dot{\b{P}}$ can thus be computed by solving another continuous time Sylvester equation involving the tangents $\b{\dot{A}}$, $\b{\dot{B}}$, and $\b{\dot{C}}$:

\begin{equation}
  \b{A}\, \b{\dot{P}} + \b{\dot{P}}\, \b{B} + 
  \big ( \b{\dot{A}} \, \b{P} + \b{P} \, \b{\dot{B}} + \b{\dot{C}} \big) = \b{0}.
\end{equation}

\paragraph{Reverse mode}
To derive the adjoints $\bar{\b{A}}, \bar{\b{B}}, \bar{\b{C}}$, we consider the Lagrangian
\begin{equation}
  \mathcal{L}(\b{A},\b{B},\b{C}, \b{P}) = \ell(\b{P}) + \text{tr}\big ( \b{S}^T (\b{A}\b{P} + \b{P}\b{B} + \b{C}) \big ),
\end{equation}
where $\ell(\b{P})$ is the scalar output w.r.t.\ which we wish to compute adjoints, and $\b{S}$ is a matrix of Lagrange multipliers used to enforce \Cref{eq:csylv}.
Taking partial derivatives on both sides of the Lagrangian and setting them to zero, we obtain:
\begin{align}
  \b{0}
   & = \frac{\partial \mathcal{L}}{\partial \b{P}}
  = \frac{\partial \ell(\b{P})}{\partial \b{P}} +
  \b{A}^T\b{S} + \b{S}\b{B}^T,\label{eq:csylv_p_reverse}          \\
  \b{0}
   & = \frac{\partial \mathcal{L}}{\partial \b{S}}
  =   \b{A}\b{P} + \b{P}\b{B} + \b{C} \label{eq:csylv_s_reverse}, \\
  \b{0}
   & = \frac{\partial \mathcal{L}}{\partial \b{A}}
  =  \b{S}\b{P}^T, \label{eq:csylv_a_reverse}                     \\
  \b{0}
   & = \frac{\partial \mathcal{L}}{\partial \b{B}}
  =  \b{P}^T\b{S},\label{eq:csylv_b_reverse}                      \\
  \b{0}
   & = \frac{\partial \mathcal{L}}{\partial \b{C}}
  =  \b{S} \label{eq:csylv_c_reverse}.
\end{align}
Note that \Cref{eq:csylv_s_reverse} does indeed enforce that $\b{P}$ be the solution of \Cref{eq:csylv}.
Identifying $\displaystyle\bar{\b{P}} = \partial \ell/\partial \b{P}$, we find that the Lagrange multiplier \(\b{S}\) satisfies another Sylvester equation
\begin{equation}
  \b{A}^T\b{S} + \b{S}\b{B}^T + \bar{\b{P}} = \b{0}.
\end{equation}
After computing $\b{S}$, and noting that when $\b{P}$ solves \Cref{eq:csylv_s_reverse} then $\mathcal{L}=\ell$, it is straightforward to compute the rest of the adjoints%
\footnote{One can formally show that $\bar{\b{X}} = \partial\mathcal{L}/\partial\b{X}$ for $\b{X} \in \{\b{A}, \b{B}, \b{C}\}$ by applying the implicit function theorem (see \citealp{Recht2016} and references therein).
}
:
\begin{equation}
  \bar{\b{A}} =\b{S} \b{P}^T,
  \qquad
  \bar{\b{B}} = \b{P}^T\b{S},
  \qquad
  \bar{\b{C}} = \b{S}.
\end{equation}
%
%Provided that the spectrum of $\b{A}$ and $\b{B}$ do not overlap and thus there is a unique solution to \Cref{eq:csylv},
%it is straightforward to verify that $\dot{P}$, $\bar{A}$, $\bar{B}$, and $\bar{C}$ are also guaranteed to be unique.

\subsection{Discrete time Sylvester equation}

We follow a similar approach as above to differentiate through the discrete time Sylvester equation, which takes the form
\begin{equation} \label{eq:dsylv}
  \b{A}\b{P}\b{B} -  \b{P} + \b{C} = \b{0},
\end{equation}
where $\b{A}$, $\b{P}$, $\b{B}$, and $\b{C}$ are all square matrices that have the same dimensions. 
There is a unique solution to \Cref{eq:dsylv} if and only if $\sigma_i(\b{A})\sigma_j(\b{B}) \neq 1$ for all $i,j$, where $\{ \sigma_i(\b{X}) \}$ denotes the eigenvalue spectrum of a matrix $\b{X}$. 

\paragraph{Forward mode}
Taking the differential on both sides of \Cref{eq:dsylv}, we have
\begin{equation}
  \b{A}\, \b{dP}\, \b{B} - \b{dP} 
  + 
  \big (
  \b{dA}\, \b{P}\b{B} 
  + 
  \b{A} \b{P} \, \b{dB} 
  + 
  \b{dC} 
  \big )
  = \b{0}.
\end{equation}
The tangent of $\b{P}$ is thus computed by solving another discrete time Sylvester equation:
\begin{equation}
  \b{A}\, \b{\dot{P}}\, \b{B} - \b{\dot{P}}
  + 
  \big (
  \b{\dot{A}}\, \b{P}\b{B} 
  + 
  \b{A} \b{P} \, \b{\dot{B}}
  + 
  \b{\dot{C}}
  \big )
  = \b{0}.
\end{equation}

\paragraph{Reverse mode}
To derive the adjoints $\bar{\b{A}}$, $\bar{\b{B}}$ and $\bar{\b{C}}$, we consider the Lagrangian
\begin{equation}
  \mathcal{L}(\b{A},\b{B},\b{C}, \b{P}, \b{S}) = \ell(\b{P}) + \text{tr}\big ( \b{S}^T (\b{A}\b{P}\b{B} - \b{P} + \b{C}) \big ).
\end{equation}
Taking partial derivatives on both sides of the Lagrangian and setting them to zero, we get:
\begin{align}
  \b{0}
   & = \frac{\partial \mathcal{L}}{\partial \b{P}}
  = \frac{\partial f(\b{P})}{\partial \b{P}} +
  \b{A}^T\b{S}\b{B}^T - \b{S},\label{eq:dsylv_p_reverse} \\
  \b{0}
   & = \frac{\partial \mathcal{L}}{\partial \b{S}}
  =   \b{A}\b{P}\b{B} - \b{P} + \b{C},                  \\
  \b{0}
   & = \frac{\partial \mathcal{L}}{\partial \b{A}}
  =  \b{S}\b{B}^T\b{P}^T, \label{eq:dsylv_a_reverse}            \\
  \b{0}
   & = \frac{\partial \mathcal{L}}{\partial \b{B}}
  =  \b{P}^T\b{A}^T\b{S},\label{eq:dsylv_b_reverse}             \\
  \b{0}
   & = \frac{\partial \mathcal{L}}{\partial \b{C}}
  =  \b{S} \label{eq:dsylv_c_reverse}.
\end{align}
Identifying $\bar{\b{P}} = \partial f(\b{P})/\partial \b{P}$, we find that the Lagrange multiplier \(\b{S}\) satisfies another discrete time Sylvester equation
\begin{equation}
  \b{A}^T\b{S}\b{B}^T - \b{S} + \bar{\b{P}} = \b{0}.
\end{equation}
The adjoints are given by
\begin{equation}
  \bar{\b{A}} =\b{S} \b{B}^T\b{P}^T,
  \quad
  \bar{\b{B}} = \b{P}^T\b{A}^T\b{S},
  \quad
  \bar{\b{C}} = \b{S}.
\end{equation}
%
%Provided that the spectrum of $\b{A}$ and $\b{B}$ do not overlap and thus there is a unique solution to \Cref{eq:dsylv},
%it is straightforward to verify that $\dot{P}$, $\bar{A}$, $\bar{B}$, and $\bar{C}$ are also guaranteed to be unique.

\section{Lyapunov equations}
\label{sec:lyap}
\subsection{Continuous time Lyapunov equation}
\label{subsec:clyap}
The continuous time Lyapunov equation is given by
\begin{equation}
  \b{A}\b{P} + \b{P}\b{A}^T + \b{Q} = \b{0}, \label{eq:clyap}
\end{equation}
where $\b{A}$, $\b{P}$, and $\b{Q}$ are square matrices that have the same dimensions.
Note that the Lyapunov equation is a special case of the Sylvester equation in \Cref{eq:csylv}, with $\b{B} = \b{A}^T$ and $\b{C} = \b{Q}$.
In control theory, we are most commonly concerned with solutions to the Lyapunov equation when $\b{Q}$ is a symmetric matrix, which also results in a symmetric $\b{P}$.
Below, we derive $\bar{\b{P}}$ and $\dot{\b{P}}$ in the general case, and specialize to the case where $\b{Q}$ is symmetric afterwards.
\paragraph{Forward mode}
Taking the differential on both sides of \Cref{eq:clyap}, we get
\begin{align}
  \b{A}\, \b{dP} + \b{dP}\, \b{A}^T +
  \big (
  \b{dA}\, \b{P} + \b{P}\, \b{dA}^T + \b{dQ}
  \big )
  = \b{0}.
\end{align}
Thus, $\b{\dot{P}}$ satisfies the continuous time Lyapunov equation:
\begin{align}
  \b{A} \b{\dot{P}} + \b{\dot{P}} \b{A}^T +
  \big (
  \b{\dot{A}}\b{P} + \b{P} \b{\dot{A}}^T + \b{\dot{Q}}
  \big )
  = \b{0}.
\end{align}
\paragraph{Reverse mode}
In reverse-mode, we consider the Lagrangian
\begin{equation}
  \mathcal{L}(\b{A},\b{Q}, \b{P}, \b{S}) = \ell(\b{P}) + \text{tr}\big (
  \b{S}^T(\b{A}\b{P} + \b{P}\b{A}^T + \b{Q}) \big ),
\end{equation}
where \(\b{S}\) is the Lagrange multiplier and $\ell(\b{P})$ is some loss downstream that is a function of $\b{P}$.
Taking partial derivatives and setting them to zero, we get the following equations
\begin{align}
  \b{0}
   & =
  \frac{\partial \mathcal{L}}{\partial \b{P}}
  =
  \frac{\partial \ell(\b{P})}{\partial \b{P}} +
  \b{A}^T\b{S} + \b{S}\b{A}             \\
  \b{0}
   & =
  \frac{\partial \mathcal{L}}{\partial \b{S}}
  =   \b{A}\b{P} + \b{P}\b{A}^T + \b{Q} \\
  \b{0}
   & =
  \frac{\partial \mathcal{L}}{\partial \b{Q}}
  =
  \b{S}                                 \\
  \b{0}
   & =
  \frac{\partial \mathcal{L}}{\partial \b{A}}
  =  \b{S}\b{P}^T + \b{S}^T\b{P}.
\end{align}
In reverse-mode, $\partial \ell(\b{P})/\partial \b{P} = \bar{\b{P}}$, which is known in the reverse pass.
We can thus compute \(\b{S}\), by solving the Lyapunov equation
\begin{equation}
  \bar{\b{P}} + \b{A}^T\b{S} + \b{S}\b{A} = \b{0}.
\end{equation}
This then gives:
\begin{equation}
  \bar{\b{A}} = \b{S}\b{P}^T + \b{S}^T \b{P}
  \quad
  \text{ and }
  \quad
  \bar{\b{Q}} = \b{S}
\end{equation}
using the equations above.
If $\b{Q}$ is a symmetric matrix, as is the case in most control applications, then $\b{P}$ and $\b{S}$ are also symmetric.
In addition, if $\b{Q}$ is positive semi-definite and all eigenvalues of $A$ have negative real parts, then $\b{P}$ and $\b{S}$ are positive semi-definite.

\subsection{Discrete time Lyapunov equation}
\label{subsec:dlyap}

The discrete time Lyapunov equation is given by
\begin{equation}
  \label{eq:dlyap}
  \b{A}\b{P}\b{A}^T - \b{P} + \b{Q} = \b{0},
\end{equation}
\paragraph{Forward mode}
Taking the differential on both sides of \Cref{eq:dlyap}, we get another discrete time Lyapunov equation:
\begin{equation}
  \b{A} \, \b{dP} \, \b{A}^T - \b{dP} +
  \big (
  \b{dA}\, \b{P}\b{A}^T + \b{A}\b{P} \,\b{dA}^T
  + \b{dQ}
  \big )
  = \b{0}.
\end{equation}

Thus, we can compute $\b{\dot{P}}$ by solving:
\begin{equation}
  \b{A} \b{\dot{P}}\b{A}^T - \b{\dot{P}} +
  \big (
  \b{\dot{A}}\b{P}\b{A}^T + \b{A}\b{P} \b{\dot{A}}^T
  + \b{\dot{Q}}
  \big )
  = \b{0}.
\end{equation}
\paragraph{Reverse mode}
In reverse-mode, we consider the following Lagrangian:
\begin{equation}
  \mathcal{L}(\b{A},\b{Q}, \b{P}, \b{S}) = \ell(\b{P}) + \text{tr}\big (\b{S}^T(\b{A}\b{P}\b{A}^T - \b{P} + \b{Q}) \big ).
\end{equation}
Taking partial derivatives and setting them to zero, we get:
\begin{align}
  \b{0}
   & =
  \frac{\partial \mathcal{L}}{\partial \b{P}}
  =
  \frac{\partial \ell(\b{P})}{\partial \b{P}} +
  \b{A}^T\b{S}\b{A} - \b{S}            \\
  \b{0}
   & =
  \frac{\partial \mathcal{L}}{\partial \b{S}}
  =   \b{A}\b{P}\b{A}^T -\b{P} + \b{Q} \\
  \b{0}
   & =
  \frac{\partial \mathcal{L}}{\partial \b{Q}}
  =
  \b{S}                                \\
  \b{0}
   & =
  \frac{\partial \mathcal{L}}{\partial \b{A}}
  =  \b{S}\b{A}\b{P}^T + \b{S}^T\b{A}\b{P}.
\end{align}
In reverse-mode, $\partial \ell(\b{P})/\partial \b{P} = \bar{\b{P}}$.
We can thus compute \(\b{S}\), by solving the discrete time Lyapunov equation
\begin{equation}
  \bar{\b{P}} + \b{A}^T\b{S}\b{A} - \b{S} = \b{0}.
\end{equation}
This then gives:
\begin{equation}
  \bar{\b{A}} = \b{S}\b{A}\b{P}^T + \b{S}^T \b{A}\b{P}
  \quad
  \text{ and }
  \quad
  \bar{\b{Q}} = \b{S}.
\end{equation}
If $\b{Q}$ is symmetric, then so are $\b{P}$ and $\b{S}$, leading to $\bar{\b{A}} = 2 \b{S}\b{P}\b{A}$.

\section{Algebraic Riccati equations}
\label{sec:are}
\subsection{Continuous time algebraic Riccati equation}
\label{subsec:care}
The continuous time algebraic Riccati equation (CARE) is given by
\begin{equation}
  \label{eq:care}
  \b{A}^T\b{P} + \b{P}\b{A} - \b{P}\b{B}\b{R}^{-1}\b{B}^T\b{P} + \b{Q} = \b{0},
\end{equation}
where $\b{Q}$ and $\b{R}$ are symmetric matrices and $\b{R}$ is positive-definite.
In most applications, the parameters of this equation will satisfy the conditions that guarantee the existence of a unique (symmetric) solution $\b{P}$.
Here, we derive $\b{\dot{P}}$ and $\b{\overline{P}}$ assuming that these conditions are met, so that $\b{P} = \b{P}^T$.

\paragraph{Forward mode}
Taking the differential on both sides of \Cref{eq:care}, we obtain a Lyapunov equation:
\begin{equation}
  \b{dP}\, \b{\tilde{A}} + \b{\tilde{A}}^T \,\b{dP}
  +\left(
  \b{dZ} + \b{dZ}^T
  + \b{dQ} + \b{K}^T\, \b{dR}\,\b{K}
  \right)
  = \b{0},
\end{equation}
where $\b{\tilde{A}} = \b{A} - \b{B}\b{K}$, $\b{K} = \b{R}^{-1}\b{B}^T\b{P}$, and $\b{dZ} = \b{P}(\b{dA} - \b{dB \, K})$.
We can thus compute $\b{\dot{\b{P}}}$ by solving the continuous time Lyapunov equation:
\begin{equation}
  \b{\dot{P}} \b{\tilde{A}} + \b{\tilde{A}}^T \b{\dot{P}}
  +\left(
  \b{\dot{Z}} + \b{\dot{Z}}^T
  + \b{\dot{Q}} + \b{K}^T\b{\dot{R}}\b{K}
  \right)
  = \b{0}
\end{equation}
where $\b{\dot{Z}} = \b{P}(\b{\dot{A}} - \b{\dot{B} \, K})$.
\paragraph{Reverse mode}
We consider the Lagrangian
\begin{equation}
  \mathcal{L}(\b{A},\b{B},\b{Q},\b{R},\b{P},\b{S}) = \ell(\b{P}) + \text{tr}\big (\b{S}^T(\b{A}^T\b{P} + \b{PA} - \b{PBR}^{-1}\b{B}^T\b{P} + \b{Q}) \big ),
\end{equation}
where $\ell(P)$ is a scalar objective function that depends on $\b{P}$ and $\b{S}$ is a Lagrange multiplier.
We are interested in computing the adjoints $\overline{\b{A}}$, $\overline{\b{B}}$, $\overline{\b{Q}}$, and $\overline{\b{R}}$.
Taking partial derivatives on both sides of the Lagrangian and setting them to zero, we get
\begin{align}
  0
   & = \frac{\partial \mathcal{L}}{\partial \b{P}}
  = \frac{\partial \ell(\b{P})}{\partial \b{P}} +
  \b{\tilde{A}} \b{S} + \b{S}\b{\tilde{A}}^T                   \\
  0
   & = \frac{\partial \mathcal{L}}{\partial \b{S}}
  =   \b{A}^T\b{P} + \b{PA} - \b{PBR}^{-1}\b{B}^T\b{P} + \b{Q} \\
  0
   & = \frac{\partial \mathcal{L}}{\partial \b{Q}}
  =  \b{S}                                                     \\
  0
   & = \frac{\partial \mathcal{L}}{\partial \b{R}}
  =  \b{K}\b{S}\b{K}^T                                         \\
  0
   & = \frac{\partial \mathcal{L}}{\partial \b{A}}
  =  \b{P} \b{S}^T + \b{P} \b{S}                               \\
  0
   & = \frac{\partial \mathcal{L}}{\partial \b{B}}
  =  -\b{PS}\b{K}^T  -\b{P}\b{S}^T\b{K}^T .
\end{align}
The Lagrange multiplier $\b{S}$ satisfies the Lyapunov equation
\begin{equation}
  \b{0} = \frac{1}{2}\left(\overline{\b{P}} + \overline{\b{P}}^T\right) + \b{\tilde{A}}\b{S} + \b{S}\b{\tilde{A}}^T \quad
  %  \text{ with }
  %  \quad
  %  \displaystyle\frac{\partial \ell(\b{P})}{\partial \b{P}} = \frac{\bar{\b{P}} + \bar{\b{P}}^T}2,
\end{equation}
where we have enforced symmetry of $\overline{\b{P}}$ because $\b{P}$ itself is symmetric -- this ensures that $\b{S}$ is also symmetric.
The adjoints are given by:
\begin{align}
  \overline{\b{A}} & = 2\b{P}\b{S}        \\
  \overline{\b{B}} & = -2\b{PS}\b{K}      \\
  \overline{\b{Q}} & = \b{S}              \\
  \overline{\b{R}} & = \b{K}\b{S}\b{K}^T.
\end{align}

\subsection{Discrete time algebraic Riccati equation}
\label{subsec:dare}
The discrete time algebraic Riccati equation (DARE) is given by
\begin{equation}
  \label{eq:dare}
  \b{A}^T\b{P}\b{A} - \b{P}
  - (\b{A}^T\b{PB})
  (\b{R} + \b{B}^T\b{P}\b{B})^{-1}(\b{A}^T \b{P}\b{B})^T + \b{Q} = \b{0}.
\end{equation}
Here, we derive $\b{\dot{P}}$ and $\b{\bar{P}}$ assuming that all conditions are satisfied that ensure the existence of a unique, symmetric solution $\b{P}$.

\paragraph{Forward mode}
Taking the differential on both sides of \Cref{eq:care} and simplifying, we get:
\begin{equation}
  \b{\tilde{A}}^T \b{dP} \, \b{\tilde{A}} - \b{dP}\\
  +
  \left(  \b{dZ} + \b{dZ}^T + \b{dQ} + \b{K}^T\,\b{dR}\,\b{K} \right)
  =
  \b{0}
\end{equation}
where $\b{\tilde{A}} = \b{A} - \b{B}\b{K}$ and $\b{K} = (\b{R} + \b{B}^T\b{PB})^{-1}\b{B}^T\b{PA}$, and $\b{dZ} = \b{\tilde{A}}^T \b{P} \left( \b{dA}- \b{dB}\,\b{K}\right)$.
Thus, in forward mode, we have the discrete time Lyapunov equation:
\begin{equation}
  \b{\tilde{A}}^T\b{\dot{P}} \b{\tilde{A}} - \b{\dot{P}}
  +
  \left( \b{\dot{Z}} + \b{\dot{Z}}^T + \b{\dot{Q}} + \b{K}^T\b{\dot{R}}\b{K} \right)
  = \b{0}.
\end{equation}
where $\dot{\b{Z}} = \b{\tilde{A}}^T \b{P} (\dot{\b{A}}- \dot{\b{B}}\b{K})$.

\paragraph{Reverse mode}
We consider the Lagrangian
\begin{equation}
  \begin{split}
    \mathcal{L}(&\b{A},\b{B},\b{Q},\b{R},\b{P},\b{S}) =\\
    &f(\b{P}) + \text{tr}\bigg (
    \b{S}^T(
    \b{A}^T\b{P}\b{A} - \b{P} -
    (\b{A}^T\b{PB})
    (\b{R} + \b{B}^T\b{P}\b{B})^{-1}(\b{A}^T\b{PB})^T + \b{Q}
    ) \bigg ).
  \end{split}
\end{equation}
Taking partial derivatives on both sides of the Lagrangian and setting them to zero, we get
\begin{align}
  0
   & = \frac{\partial \mathcal{L}}{\partial \b{P}}
  = \frac{\partial f(\b{P})}{\partial \b{P}} +
  \b{\tilde{A}} \b{S} \b{\tilde{A}}^T - \b{S}              \\
  0
   & = \frac{\partial \mathcal{L}}{\partial \b{S}}
  =
  \b{A}^T\b{P}\b{A} - \b{P} -
  (\b{A}^T\b{PB})
  (\b{R} + \b{B}^T\b{P}\b{B})^{-1}(\b{A}^T \b{PB})^T + \b{Q}  \\
  0
   & = \frac{\partial \mathcal{L}}{\partial \b{Q}}
  =  \b{S}                                                 \\
  0
   & = \frac{\partial \mathcal{L}}{\partial \b{R}}
  =  \b{K}\b{S}\b{K}^T
  \\
  0
   & = \frac{\partial \mathcal{L}}{\partial \b{A}}
  =  \b{P} \b{\tilde{A}} \b{S}^T + \b{P}\b{\tilde{A}}\b{S} \\
  0
   & = \frac{\partial \mathcal{L}}{\partial \b{B}}
  =  -\b{P\tilde{A}SK}^T  -\b{P\tilde{A}}\b{S}^T\b{K}^T
\end{align}
The Lagrange multiplier $\b{S}$ satisfies the Lyapunov equation
\begin{equation}
  \b{0} = \frac{1}{2}(\bar{\b{P}} + \bar{\b{P}}^T) + \b{\tilde{A}}\b{S}\b{\tilde{A}}^T - \b{S}
%  \text{ with }
%  \quad
%  \displaystyle\frac{\partial f(\b{P})}{\partial \b{P}} = \frac{1}2 (\bar{\b{P}} + \bar{\b{P}}^T),
\end{equation}
where we have enforced the symmetry of $\bar{\b{P}}$ because $\b{P}$ is symmetric.
This ensures that $\b{S}$ is also symmetric and the adjoints are given by:
\begin{align}
  \bar{\b{A}} & = 2\b{P}\b{\tilde{A}}\b{S} \\
  \bar{\b{B}} & = -2\b{P\tilde{A}SK}^T     \\
  \bar{\b{Q}} & = \b{S}                    \\
  \bar{\b{R}} & = \b{K}\b{S}\b{K}^T.
\end{align}

\section{Example application: inverse LQR}
\label{sec:examples}
\begin{figure}[t]
    \centering
    \includegraphics[width=0.75\textwidth]{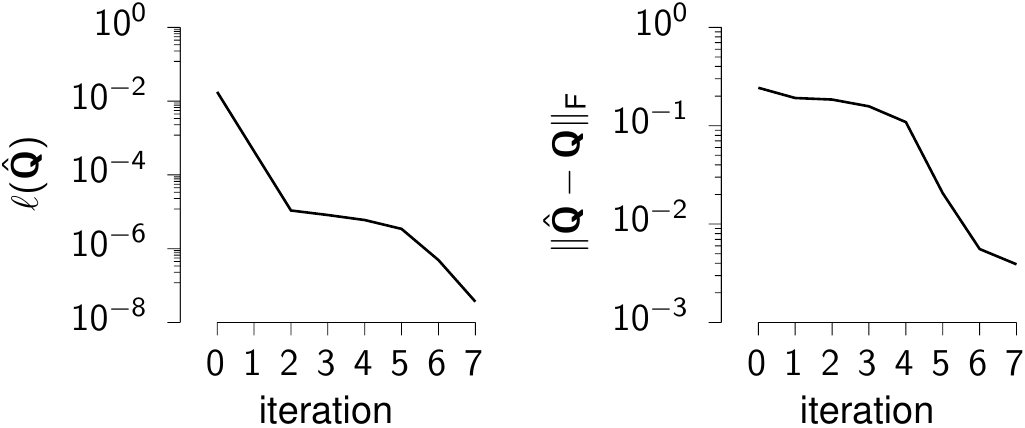}
    \caption{\label{fig:example}
    Evolution of $\ell(\b{\hat{Q}})$ (left) and $\|\b{\hat{Q}}\|$ (right) duing optimization (see text).
    }
\end{figure}

To show a concrete application of automatic differentiation through these algebraic matrix equations, we consider a discrete time linear system that evolves according to the equations
\begin{equation}
    \b{x}_{t+1} = \b{A}\b{x}_t + \b{B}\b{u}_t,
    \quad
    \b{x}_0 = \b{x}_{\text{init}}
\end{equation}
where $\b{x}_t$ is the state of the system, $\b{u}_t$ is
some input that enters the system through some matrix $\b{B}$, and $\b{A}$ is the state transition matrix.
The infinite-horizon, discrete-time LQR problem involves finding the optimal inputs $\{ \b{u}_0, \b{u}_1, \ldots \}$ that minimise \citep{Anderson2007} 
\begin{equation}
    \label{eq:lqr_cost}
    J[\b{u}_0, \b{u}_1, \cdots] =
    \sum_{t=0}^\infty
    \b{x}_t^T \b{Q} \b{x}_t + \b{u}_t^T \b{R} \b{u}_t
\end{equation}
for some symmetric positive semi-definite matrices $\b{Q}$ and $\b{R}$.
It is well-known that the optimal solution is linear state-feedback:
\begin{equation}
    \b{u}_t = -\b{K}\b{x}_t,
\end{equation}
where $\b{K} = (\b{R} + \b{B}^T\b{PB})^{-1}\b{B}^T\b{PA}$ and $\b{P}$ satisfies a DARE (c.f.\ \Cref{eq:dare}):
\begin{equation}
    \b{A}^T\b{P} + \b{P}\b{A} - \b{PBR}^{-1}\b{B}^T\b{P} + \b{Q} = \b{0}.
\end{equation}
We consider the following inverse LQR problem: can we infer $\b{Q}$, i.e.\ the way state deviations from zero are penalized in various state space directions, given sample state trajectories generated by the optimally controlled system, and given knowledge of $\b{A}$, $\b{B}$, and $\b{R}$?
More formally, given $K$ trajectories of the system under optimal LQR control sampled at $T$ time points:
\begin{equation}
    \b{x}^{(k)}_{0}, \b{x}^{(k)}_{1}, \cdots \b{x}^{(k)}_{T-1} \quad \text{ for } \quad k = 1, \cdots, K,
\end{equation}
can we determine the matrix $\b{Q}$ used in the LQR algorithm's objective function (\Cref{eq:lqr_cost})?
Assuming that this problem is well-posed (i.e., $\b{Q}$ is identifiable; \citealp{Zhang2019}), we can solve this problem by differentiating through DARE (\Cref{subsec:dare}).
More specifically, we start with an initial guess of $\b{Q}$, which we denote as $\b{\hat{Q}}$.
Next, we solve the LQR problem and find the corresponding optimal inputs $\b{\hat{u}}_t$ and trajectories $\b{\hat{x}}_t$.
We then minimize the objective
\begin{equation}
    \mathcal{\ell}(\b{\hat{Q}})=
    \frac{1}{KT} \sum_{k=1}^{K} \sum_{t=0}^{T-1}
    (\b{x}^{(k)}_{t} - \b{\hat{x}}^{(k)}_{t})^2.
\end{equation}
with respect to $\b{\hat{\b{Q}}}$.
To illustrate our solution, we consider the following example system:
\begin{equation}
    \begin{split}
        \b{A} =
        \begin{pmatrix}
            1 & 1 \\
            0 & 1
        \end{pmatrix}
        \quad
        &
        \quad
        \b{B} =
        \begin{pmatrix}
            1 & 0\\
            0 & 1
        \end{pmatrix}\\
        \b{Q} =
        \begin{pmatrix}
            1 & 0 \\
            0 & 0
        \end{pmatrix}
        \quad
        &
        \quad
        \b{R} =
        \begin{pmatrix}
            0.1 & 0   \\
            0   & 0.3
        \end{pmatrix}.
    \end{split}
\end{equation}
We created a synthetic data set by sampling $\b{K}=30$ initial states $\b{x}_{\text{init}}$ from $\mathcal{N}(\b{0}, \b{I})$ and for each sample, we simulated the activity trajectories under optimal control for $T=30$ steps.
We then optimized $\ell(\b{\hat{Q}})$ for the synthetic data set using L-BFGS~\citep{Liu89}.
After $7$ iterations, we were able to recover the true $\b{Q}$ to a good degree of accuracy (\Cref{fig:example}).
Code for this example is available at \href{https://github.com/tachukao/autodiff-inverse-lqr}{https://github.com/tachukao/autodiff-inverse-lqr}.

\section{Conclusion}
\label{sec:conclusion}
In this note, we derived the forward and reverse mode derivatives for the solutions to the discrete time and continuous time variants of the Sylvester, Lyapunov, and Riccati equations (summarized in \Cref{sm:table:forward,sm:table:reverse}).
These equations are widely-used in control theory and other branches of applied mathematics.
We demonstrate the usefulness of these derivatives on an inverse LQR problem, where parameters of the optimal control problem are inferred from observations of state trajectories by differentiating through the discrete time Riccati equation (\Cref{subsec:dare}).
These derivatives (with the exception of the discrete time Sylvester equation, of more limited use), have been implemented and numerically tested in Owl, a numerical library written in OCaml with a full-featured automatic differentiation module.

\section*{Acknowledgements}
\label{sec:acknowledgements}
We are grateful for helpful comments on the manuscript by Kris Jensen.

\bibliographystyle{apalike}
\bibliography{references}

\newpage
\begin{appendix}
  \section{Forward-mode derivatives summary}
  % Forward-mode Summary
  \begin{table}[!h]
  \small
  \centering \rowcolors{1}{gray!5}{white}
  \setlength\extrarowheight{3pt}
  \begin{tabular}{|c|c|}
    \rowcolor{brown!10}
    \hline
    \textbf{Equation} &
    \textbf{Tangents}
    \\
    %\noalign{\bigskip}
    \multicolumn{2}{c}{
      \textsc{\textcolor{black}{
          continuous time Sylvester
        }}
      \Cref{eq:csylv}
    }
    \\[0.5em]
    $\scriptstyle \b{A}\b{P}+ \b{P}\b{B} + \b{C} = \b{0}$
                      &
    $\scriptstyle \b{A}\, \b{\dot{P}} + \b{\dot{P}}\, \b{B} +
      \big ( \b{\dot{A}} \, \b{P} + \b{P} \, \b{\dot{B}} + \b{\dot{C}} \big) = \b{0}$
    \\
    %
    %\noalign{\bigskip}
    \multicolumn{2}{c}{
      \textsc{\textcolor{black}{
          discrete time Sylvester
        }}
      \Cref{eq:dsylv}
    }
    \\[0.5em]
    $\scriptstyle \b{A}\b{P}\b{B} -  \b{P} + \b{C} = \b{0}$
                      &
    $\scriptstyle
      \b{A}\, \b{\dot{P}}\, \b{B} - \b{\dot{P}}
      +
      \big (
      \b{\dot{A}}\, \b{P}\b{B}
      +
      \b{A} \b{P} \, \b{\dot{B}}
      +
      \b{\dot{C}}
      \big )
      = \b{0}$
    \\
    %
    %\noalign{\bigskip}
    \multicolumn{2}{c}{
      \textsc{\textcolor{black}{
          continuous time Lyapunov
        }}
      \Cref{eq:clyap}
    }
    \\[0.5em]
    $\scriptstyle
      \b{A}\b{P}+ \b{P}\b{A}^T + \b{C} = \b{0}$
                      &
    $\scriptstyle
      \b{A}\, \b{\dot{P}} + \b{\dot{P}}\, \b{A}^T +
      \big ( \b{\dot{A}} \, \b{P} + \b{P} \, \b{\dot{A}}^T + \b{\dot{C}} \big) = \b{0}$
    \\
    %
    %\noalign{\bigskip}
    \multicolumn{2}{c}{
      \textsc{\textcolor{black}{
          discrete time Lyapunov
        }}
      \Cref{eq:dlyap}
    }
    \\[0.5em]
    $\scriptstyle
      \b{A}\b{P}\b{A}^T -  \b{P} + \b{C} = \b{0}$
                      &
    $
      \scriptstyle
      \b{A}\, \b{\dot{P}}\, \b{A}^T - \b{\dot{P}}
      +
      \big (
      \b{\dot{A}}\, \b{P}\b{A}^T
      +
      \b{A} \b{P} \, \b{\dot{A}}^T
      +
      \b{\dot{C}}
      \big )
      = \b{0}$
    \\
    %\noalign{\bigskip}
    \multicolumn{2}{c}{
      \textsc{\textcolor{black}{
          continuous time algebraic Riccati
        }}
      \Cref{eq:care}
    }
    \\[0.5em]
    $
      \scriptstyle
      \b{A}^T\b{P} + \b{P}\b{A} - \b{P}\b{B}\b{R}^{-1}\b{B}^T\b{P} + \b{Q} = \b{0}$
                      &
    $
      \begin{array}{rl}
        \scriptstyle
        \b{\dot{P}} \b{\tilde{A}} + \b{\tilde{A}}^T \b{\dot{P}}
        +\big (
        \b{\dot{A}}^T \b{P} + \b{P} \b{\dot{A}}
        - \b{P\dot{B}K}
        - \b{K}^T\b{\dot{B}}^T\b{P}
         & \\
        \scriptstyle
        + \b{\dot{Q}} + \b{K}^T\b{\dot{R}}\b{K}
        \big)
        =\b{0}
         &
      \end{array}$
    \\
    %
    %\noalign{\bigskip}
    \multicolumn{2}{c}{
      \textsc{\textcolor{black}{
          discrete time algebriac Riccati
        }}
      \Cref{eq:dare}
    }
    \\[0.5em]
    $
      \begin{array}{rr}
        \scriptstyle
        \b{A}^T\b{P}\b{A} - \b{P} & \\
        \scriptstyle
        - (\b{A}^T\b{PB})(\b{R} + \b{B}^T\b{P}\b{B})^{-1}(\b{B}^T\b{PA})
        + \b{Q}
        = \scriptstyle
        \b{0}
      \end{array}$
                      &

    $
      \begin{array}{rl}
        \scriptstyle
        \b{\tilde{A}}^T\b{\dot{P}} \b{\tilde{A}} - \b{\dot{P}}
        +
        \big (
        \b{\dot{A}}^T\b{P}\b{\tilde{A}}
        + \b{\tilde{A}}^T\b{P}\b{\dot{A}}
         & \\
        \scriptstyle
        - \b{K}^T\b{\dot{B}}^T\b{P\tilde{A}}
        - \b{\tilde{A}}^T\b{P}\b{\dot{B}}\b{K}
        + \b{\dot{Q}} + \b{K}^T\b{\dot{R}}\b{K}
        \big )
        = \b{0}
      \end{array}$
    \\
    \hiderowcolors
    \noalign{\bigskip}
    \hline
  \end{tabular}
  \caption{\label{sm:table:forward}Summary of forward-mode derivatives}
\end{table}
  \newpage
  \section{Reverse-mode derivatives summary}
  % Reverse-mode Summary
  \begin{table}[!h]
    \small
    \centering \rowcolors{1}{gray!5}{white}
    \setlength\extrarowheight{3pt}
    \begin{tabular}{|c|c|}
        \rowcolor{brown!10}
        \hline
        \textbf{Equation} &
        \textbf{Adjoints}
        \\
        %\noalign{\bigskip}
        \multicolumn{2}{c}{
            \textsc{\textcolor{black}{
                    continuous time Sylvester
                }}
            \Cref{eq:csylv}
        }
        \\[0.5em]
        $\scriptstyle \b{A}\b{P}+ \b{P}\b{B} + \b{C} = \b{0}$
                          &

        $
            \begin{array}{rl}
                 & \scriptstyle
                \b{0}=
                \b{A}^T\b{S} + \b{S}\b{B}^T + \bar{\b{P}} \\
                 &
                \scriptstyle
                \bar{\b{A}}
                = \b{S} \b{P}^T                           \\
                 &
                \scriptstyle
                \bar{\b{B}}
                = \b{P}^T\b{S}                            \\
                 &
                \scriptstyle
                \bar{\b{C}}
                = \b{S}
            \end{array}
        $
        \\
        %
        %\noalign{\bigskip}
        \multicolumn{2}{c}{
            \textsc{\textcolor{black}{
                    discrete time Sylvester
                }}
            \Cref{eq:dsylv}
        }
        \\[0.5em]
        $\scriptstyle \b{A}\b{P}\b{B} -  \b{P} + \b{C} = \b{0}$
                          &
        $
            \begin{array}{rl}
                 &
                \scriptstyle
                \b{0} = \b{A}^T\b{S}\b{B}^T - \b{S} + \bar{\b{P}}, \\
                 & \scriptstyle
                \bar{\b{A}} =\b{S} \b{B}^T\b{P}^T                  \\
                 & \scriptstyle
                \bar{\b{B}} = \b{P}^T\b{A}^T\b{S}                  \\
                 & \scriptstyle
                \bar{\b{C}} = \b{S}
            \end{array}
        $
        \\
        %
        %\noalign{\bigskip}
        \multicolumn{2}{c}{
            \textsc{\textcolor{black}{
                    continuous time Lyapunov
                }}
            \Cref{eq:clyap}
        }
        \\[0.5em]
        $\scriptstyle
            \b{A}\b{P}+ \b{P}\b{A}^T + \b{C} = \b{0}$
                          &

        $
            \begin{array}{rl}
                 &
                \scriptstyle
                \b{0} =
                \bar{\b{P}} + \b{A}^T\b{S} + \b{S}\b{A}    \\
                 & \scriptstyle
                \bar{\b{A}} = \b{S}\b{P}^T + \b{S}^T \b{P} \\
                 & \scriptstyle
                \bar{\b{Q}} = \b{S}
            \end{array}
        $

        \\
        %
        %\noalign{\bigskip}
        \multicolumn{2}{c}{
            \textsc{\textcolor{black}{
                    discrete time Lyapunov
                }}
            \Cref{eq:dlyap}
        }
        \\[0.5em]
        $\scriptstyle
            \b{A}\b{P}\b{A}^T -  \b{P} + \b{C} = \b{0}$
                          &
        $
            \begin{array}{rl}
                 & \scriptstyle
                \b{0}=
                \bar{\b{P}} + \b{A}^T\b{S}\b{A} - \b{S}              \\
                 & \scriptstyle
                \bar{\b{A}} = \b{S}\b{A}\b{P}^T + \b{S}^T \b{A}\b{P} \\
                 & \scriptstyle
                \bar{\b{Q}} = \b{S}.
            \end{array}
        $
        \\
        %\noalign{\bigskip}
        \multicolumn{2}{c}{
            \textsc{\textcolor{black}{
                    continuous time algebraic Riccati
                }}
            \Cref{eq:care}
        }
        \\[0.5em]
        $
            \scriptstyle
            \b{A}^T\b{P} + \b{P}\b{A} - \b{P}\b{B}\b{R}^{-1}\b{B}^T\b{P} + \b{Q} = \b{0}$
                          &
        $
            \begin{array}{rl}

                 & \scriptstyle
                \b{0} = \frac{1}{2}(\bar{\b{P}} + \bar{\b{P}}^T) + \b{\tilde{A}}\b{S} + \b{S}\b{\tilde{A}}^T \\
                 & \scriptstyle
                \bar{\b{A}} = 2\b{P}\b{S}                                                                    \\
                 & \scriptstyle
                \bar{\b{B}} = -2\b{PS}\b{K}                                                                  \\
                 & \scriptstyle
                \bar{\b{Q}} = \b{S}                                                                          \\
                 & \scriptstyle
                \bar{\b{R}} = \b{K}\b{S}\b{K}^T.
            \end{array}$
        \\
        %
        %\noalign{\bigskip}
        \multicolumn{2}{c}{
            \textsc{\textcolor{black}{
                    discrete time algebriac Riccati
                }}
            \Cref{eq:dare}
        }
        \\[0.5em]
        $
            \begin{array}{rr}
                \scriptstyle
                \b{A}^T\b{P}\b{A} - \b{P} & \\
                \scriptstyle
                - (\b{A}^T\b{PB})(\b{R} + \b{B}^T\b{P}\b{B})^{-1}(\b{B}^T\b{PA})
                + \b{Q}
                = \scriptstyle
                \b{0}
            \end{array}$
                          &

        $
            \begin{array}{rl}
                 & \scriptstyle \b{0} = \frac{1}{2}(\bar{\b{P}} + \bar{\b{P}}^T) + \b{\tilde{A}}\b{S}\b{\tilde{A}}^T - \b{S} \\
                 & \scriptstyle
                \bar{\b{A}} = 2\b{P}\b{\tilde{A}}\b{S}                                                                \\
                 & \scriptstyle
                \bar{\b{B}} = -2\b{P\tilde{A}SK}^T                                                                    \\
                 & \scriptstyle
                \bar{\b{Q}} = \b{S}                                                                                   \\
                 & \scriptstyle
                \bar{\b{R}} = \b{K}\b{S}\b{K}^T
            \end{array}$
        \\
        \hiderowcolors
        \hline
    \end{tabular}
    \caption{\label{sm:table:reverse}Summary of reverse-mode derivatives}
\end{table}
\end{appendix}

\end{document}